\documentclass[11pt,leqno]{amsart}

\usepackage[top=2.5 cm,bottom=2 cm,left=2.5 cm,right=2 cm]{geometry}
\usepackage[utf8]{inputenc}
\usepackage{xfrac}
\usepackage{bbm}
\usepackage{graphicx}
\usepackage{hyperref}
\hypersetup{
    hidelinks,
    }
\usepackage{color}
\usepackage{amssymb}
\usepackage{amsfonts}
\usepackage{psfrag}
\newcounter{stepnb}

\usepackage{tikz}
\usetikzlibrary{shapes.geometric,calc,decorations.markings,math}

\tikzstyle{nodo}=[circle,draw,fill,inner sep=0pt,minimum size=%
1.5mm]
\tikzstyle{infinito}=[circle,inner sep=0pt,minimum size=0mm]

\newtheorem{theorem}{Theorem}[section]
\newtheorem{lemma}[theorem]{Lemma}

\newtheorem{definition}[theorem]{Definition}
\numberwithin{equation}{section}

\newcommand{\R}{\mathbb{R}}

\newcommand{\be}{\begin{equation}}
\newcommand{\eq}{\end{equation}}

\begin{document}
\title[]{Well-posedness results for the Generalized Aw-Rascle-Zhang model}

\author[E. Marconi]{Elio Marconi}
\address{E.M. Dipartimento di Matematica ``Tullio Levi Civita", Universit\`a degli Studi di Padova, Via
Trieste 63, 35131 Padova, Italy}
\email{elio.marconi@unipd.it}
\author[L.~V.~Spinolo]{Laura V.~Spinolo}
\address{L.V.S. CNR-IMATI ``E. Magenes", via Ferrata 5, I-27100 Pavia, Italy.}
\email{spinolo@imati.cnr.it}
\maketitle
{
\rightskip .85 cm
\leftskip .85 cm
\parindent 0 pt
\begin{footnotesize}
We establish existence, uniqueness and stability results for the so-called Generalized Aw-Rascle-Zhang model, a second order traffic model introduced by Fan, Herty and Seibold. Our analysis is motivated by the companion paper~\cite{MS:nonlocalGARZ}. 

\medskip



\end{footnotesize}
}

\section{Introduction}
This note deals with the Generalized Aw-Rascle-Zhang model introduced in~\cite{FanHertySeibold}, namely 
\be \label{e:GARZ}
\left\{
\begin{array}{ll}
        \partial_t \rho + \partial_x [V(\rho, u) \rho] =0 \\
        \partial_t u + V(\rho, u) \partial_x u = 0. \\
\end{array}
\right.
\eq
In the previous expression, the unknown $\rho: \R_+ \times \R \to \R$ denotes the density of cars, and $V$ their velocity, whereas the unknown $u: \R_+ \times \R \to \R$ represents the driving style of drivers (for instance, their empty road velocity, that is the velocity they choose when the road is completely free). As such, it is a Lagrangian marker governed by the transport equation at the second line of~\eqref{e:GARZ}, whereas the first line of~\eqref{e:GARZ} expresses the conservation of the total amount of cars. System~\eqref{e:GARZ} is formally equivalent to the system of conservation laws
\be \label{e:GARZcf}
\left\{
\begin{array}{ll}
        \partial_t \rho + \partial_x [V(\rho, u) \rho] =0 \\
        \partial_t [\rho u] + \partial_x [V(\rho, u) \rho u] = 0, \\
\end{array}
\right.
\eq
which, contrary to~\eqref{e:GARZ}, admits a standard notion of solution in the sense of distributions. In the following we will be only concerned with solutions $(\rho, u)$ such that $u$ is a Lipschitz continuous function, and hence we will use the formulation~\eqref{e:GARZ} rather than~\eqref{e:GARZcf}: the second equation of~\eqref{e:GARZ} will be satisfied as an identity between $L^\infty$ functions. The use of~\eqref{e:GARZ} eases the analysis of the vanishing $\rho$ case, which requires special care when one uses~\eqref{e:GARZcf}. 

In this note we deal with the Cauchy problem posed by coupling~\eqref{e:GARZ} with the initial data
\be \label{e:idpose}
     \rho(0, \cdot) = \rho_0, \qquad u(0, \cdot) = u_0,
\eq
and, as in~\cite{FanHertySeibold}, we assume, in view of modeling considerations, that 
\begin{equation} \label{e:id}
0 \leq \rho_0 \leq 1, \quad u_0 \in L^\infty (\R), \; u_0 \ge 0
\end{equation}
and that\footnote{To simplify the exposition we require that the conditions~\eqref{e:V} are satisfied on the whole $\R^2$. In the following however we show that the entropy admissible solution $(\rho, u) \in [0, 1]\times [0, \| u_0 \|_{L^\infty}]$ and hence as a matter of fact it suffices to require that the conditions are satisfied in that range.}
\be \label{e:V}
         V\in C^2 (\R^2; \R), \quad V \ge 0, \quad 
    \partial_1 V \leq 0, \quad  \partial_2 V \ge 0, \quad V(1, w) =0 \; \text{for every $w$},
\eq 
where $\partial_1V$ and $\partial_2V$ denote the partial derivatives of $V$ with respect to $\rho$ and $u$, respectively. In the previous equations and in the following we normalize to $1$ the maximal possible car density, corresponding to bumper-to-bumper packing. 

In the present note we establish well-posedness results for~\eqref{e:GARZ},\eqref{e:idpose} in a suitable class of functions. Our analysis is motivated by the companion paper~\cite{MS:nonlocalGARZ}, to which we refer for a wider description of the model and a more complete list of references. Here we only refer to~\cite{GaravelloHanPiccoli} for an in-depth introduction to traffic flow models. 

First, we define our notion of solution of~\eqref{e:GARZ}.
\begin{definition} \label{d}
We term $(\rho, u) \in L^\infty_{\mathrm{loc}}  (\R_+\times \R) \cap C^0 ( \R_+, L^1_{\mathrm{loc}} (\R)) \times W^{1 \infty} (\R_+, \R))$ an \emph{entropy admissible solution} of~\eqref{e:GARZ} if the equation at the first line is satisfied in the sense of distributions, the second as an equality between $L^\infty$ functions, and the following entropy condition holds: 
\be 
\label{e:entropy}
   \partial_t |\rho- k| + \partial_x \Big[\mathrm{sign}[\rho-k] \big[V(\rho, u) \rho - V(k, u) k \big]\Big] +
    \mathrm{sign}[\rho-k] k \partial_2 V (k, u) \partial_x u \leq 0
\eq
in the sense of distributions on $\R_+ \times \R$, for every $k \in \R$. 
\end{definition}
Note that, since $\rho \in C^0 ( \R_+, L^1_{\mathrm{loc}} (\R))$ and $u$ is a Lipschitz continuous function, we can give a meaning to the values $\rho(t, \cdot)$ and $u(t, \cdot)$ at \emph{every} $t$. An entropy admissible solution of the Cauchy problem~\eqref{e:GARZ},\eqref{e:idpose} is then an entropy admissible solution of~\eqref{e:GARZ} which attains the initial datum~\eqref{e:idpose}.

With the above definition in place, we can now state our well-posedness result. 
\begin{theorem}
\label{t:wpl}
Assume $V$ satisfies~\eqref{e:V} and that the initial data $(\rho_0, u_0)$ satisfy~\eqref{e:id}, $\rho_0 \in BV (\R)$ and 
\be \label{e:id2}
    u_0 \in W^{1 \infty} (\R), \quad  u_0' = z_0 \rho_0, \quad z_0 \in L^\infty (\R)
\eq
and
\be \label{e:idgarz}
      z_0' = \rho_0 \psi_0, \quad \text{for some $\psi_0 \in L^\infty (\R)$}. 
\eq
Then there is a unique entropy admissible solution of \eqref{e:GARZ},\eqref{e:idpose} in the class of functions $(\rho, u) \in L^\infty_{\mathrm{loc}} (\R_+;BV(\R)) \times  L^\infty (\R_+; W^{1 \infty}(\R))$ such that $\partial_x u = \rho z$ with $z \in W^{1 \infty} (\R_+ \times \R)$.  Also, this unique solution satisfies 
\begin{eqnarray}
    0 \leq \rho \leq 1, \; \mathrm{TotVar} \, \rho(t, \cdot) \leq D(t), \; \text{for every $t \ge 0$},  \\
 \label{e:verde3} 0 \leq u_0 \leq \| u_0 \|_{L^\infty}, \quad    \| z \|_{L^\infty} \leq \| z_0 \|_{L^\infty}, \quad  \partial_x z = \rho \psi, \quad \| \psi \|_{L^\infty} \leq \| \psi_0 \|_{L^\infty},
 \end{eqnarray}
 where $D(t)$ is a suitable constant depending on $t,  V, \| \rho_0 \|_{L^1}, \| z_0 \|_{L^\infty}, \| \psi_0 \|_{L^\infty}$ and $\mathrm{TotVar}\, \rho_0$. Given $T>0$, we also have the stability estimate 
\be \label{e:stability}
      \| u_1 (t, \cdot) - u_2 (t, \cdot) \|_{C^0} + \| \rho_1(t, \cdot) - \rho_2 (t, \cdot) \|_{L^1}
      \leq K \Big[
       \| u_1 (0, \cdot) -  u_2 (0, \cdot) \|_{C^0} + \| \rho_1 (0, \cdot)- \rho_2 (0, \cdot)
     \|_{L^1} \Big]
\eq
for every $t \in [0, T]$, where $K$ is a suitable constant only depending on $V$ and on the quantities $T,  \| \rho_i(0, \cdot) \|_{L^1}$, $ \| z_i (0, \cdot) \|_{L^\infty}, \| \psi_i (0, \cdot) \|_{L^\infty}$, $\mathrm{TotVar}\, \rho_i (0, \cdot)$ for $i=1, 2$.
\end{theorem}
As mentioned before, the primary motivation for our analysis comes from~\cite{MS:nonlocalGARZ}, where we recover~\eqref{e:GARZ} as the local limit of a nonlocal model. In the proof of the uniqueness part we adapt an argument due to Tveito and Winther~\cite{TveitoWinther}, whereas the proof
of the existence part is based on an iteration argument as in~\cite{MS:nonlocalGARZ}. We conclude this introduction with a comparison with some previous works. 
\subsection*{Comparison with works on a model for polymer flooding in oil recovery}
System~\eqref{e:GARZ} has been extensively studied as a model for a multiphase flow used in oil recovery, see for instance~\cite{IsaacsonTemple,Temple_oil,TveitoWinther} and~\cite{Shen}, the last one also containing and in-depth discussion on an extended list of references.

As pointed ou for instance in~\cite[\S3]{Temple_oil}, \eqref{e:GARZ} is an hyperbolic system of conservation laws in the variables $(\rho, m = \rho u)$. The eigenvalues of the Jacobian 
matrix of the flux are $\lambda_1 = V+ \rho \partial_1 V$ and $\lambda_2 = V$, and hence the systems is strictly hyperbolic if for instance $\rho$ is positive and bounded away from and $0$ and $\partial_1 V <0$, a strictly stronger condition than the third one in~\eqref{e:V}. The second characteristic field belongs to the Temple class (that is, integral curves of the eigenvectors are straight lines), and is genuinely nonlinear under assumptions that are fairly reasonable in the traffic flow framework. The second vector field is linearly degenerate. When the system is strictly hyperbolic, the by-now classical Glimm-Bressan theory applies, see~\cite{Bressan_book}, and yields global-in-time existence and uniqueness for small total variation data.  

By using the special structure of~\eqref{e:GARZ}, various works like~\cite{IsaacsonTemple,Temple_oil} established global-in-time existence results under much more general assumptions on the data. The reason why we could not directly apply these results is because they impose conditions on the function $V$ that are natural in view of applications to oil recovery, but not so much in the vehicular traffic framework, for instance they typically require that the function $\rho \mapsto \rho V(\rho, u)$ is s-shaped.  Note, however, that these assumptions imply that there is a curve inside the range of values that the solutions attain where the eigenvalues switch order (that is, $\lambda_1 < \lambda_2$ on one side of the curve, and $\lambda_1 > \lambda_2$ on the other side). Conversely, in the traffic flow framework, strictly hyperbolicity is lost on the boundary of the phase space domain, which in principle should make the analysis easier. It seems, therefore, reasonable to expect that one could, for instance, adapt the techniques of~\cite{IsaacsonTemple,Shen,Temple_oil}, and obtain global-in-time existence results under assumptions that are natural in view of applications to traffic flows. 
Note, however, that the uniqueness result in~\cite{TveitoWinther} requires much stronger regularity on $u$ than the one given by the existence results in~\cite{IsaacsonTemple,Temple_oil}, and this higher regularity is also what we need in the nonlocal-to-local limit result in~\cite{MS:nonlocalGARZ}. We also refer to~\cite{Shen_uniqueness} for uniqueness results for the Riemann problem. 

Wrapping up, the novelties of the present note are that we adapt the uniqueness result in~\cite{TveitoWinther} to cover assumptions that are natural in the framework of traffic models, and that we provide a new global-in-time existence proof (under restrictive assumptions on the initial data), based on a argument very different from the one in~\cite{Temple_oil}. Using the formulation~\eqref{e:GARZ} allows us to avoid the problems stemming from the fact that the function $u$ obtained using~\eqref{e:GARZcf} is not uniquely defined on the sets where $\rho$ vanishes. We also refer to the ongoing project~\cite{ACCGK} for other results on the Generalized Aw-Rascle-Zhang system. 
\subsection*{Outline} The exposition is organized as follows. In~\S\ref{s:ex} we establish the existence part of Theorem~\ref{t:wpl}, in \S\ref{ss:luni2} the uniqueness and stability part. 
\subsection*{Notation} We denote by $C(a_1, \dots, a_n)$ a constant only depending on the quantities $a_1, \dots, a_n$. Its precise value can vary from occurrence to occurrence. 
\section{Existence}\label{s:ex}
To establish the existence part of Theorem~\ref{t:wpl} proceed as follows: in \S\ref{sss:lapp} we construct the approximating sequence, in \S\ref{sss:lpass} we pass to the limit and in \S\ref{ss:pphienne} we establish the proof of Lemma~\ref{l:phienne}.
\subsection{Construction of the approximation}\label{sss:lapp}
First of all, we point out that with no loss of regularity we can assume that the initial datum is compactly supported, that is 
\be \label{e:tortora}
    \rho_0 (x) = 0 \; \text{for a.e $x \notin [-R, R]$ for some $R>0$.} 
\eq
Indeed, once we establish the existence part of Theorem~\ref{t:wpl} under the further assumption~\eqref{e:tortora}, we can remove it by relying on a fairly standard approximation argument.  We then proceed according to the following steps. \\
{\bf Step 1:} we set $\rho_1(t, x) : = \rho_0(x)$, $u_1(t, x)= u_0(x)$. Next, we fix $\tau_0$, to be determined in the following, and $(\rho_n, u_n) $ satisfying 
\be \label{e:hpitera1}
    0 \leq \rho_n \leq 1 \; \text{on $]0, \tau_0[ \times \R$}, \quad  \| \rho (t, \cdot) \|_{L^1 (\R)} \leq \| \rho_0 \|_{L^1}, \quad 
       \mathrm{TotVar} \rho_{n}(t, \cdot) \leq M_0 \quad \text{for every $t \in ]0, \tau_0[$},
\eq 
with $M_0$ to be determined in the following, 
and 
\begin{eqnarray} 
    \| u_n \|_{L^\infty} \leq \| u_0 \|_{L^\infty}, \quad \partial_x u_n = \rho_n z_n, 
    \quad \| z_n \|_{L^\infty} \leq \| z_0 \|_{L^\infty} \label{e:hpitera2}\\
          \partial_x z_n = \rho_n \psi_n, \quad \| \psi_n \|_{L^\infty} \leq  \| \psi_0 \|_{L^\infty} \label{e:hpitera3} .
\end{eqnarray}
On the time interval $t \in ]0, \tau_0[$, we construct $\rho_{n+1}$ as the entropy admissible solution of the Cauchy problem
\be \label{e:kruenne}
    \left\{
   \begin{array}{ll}
              \partial_t \rho_{n+1} + \partial_x [V(\rho_{n+1}, u_n) \rho_{n+1}] =0 \\
              \rho_{n+1}(0, \cdot) = \rho_0, \\
   \end{array}
   \right.
\eq
which immediately yields 
\be \label{e:rhon+1elle1}
     \| \rho_{n+1} (t, \cdot) \|_{L^1} \leq \| \rho_0 \|_{L^1}.
\eq
{\bf Step 2:} we show that 
\be \label{e:bdfrombelow}
     0  \leq \rho_{n+1}  \leq 1 \; \text{on $]0, \tau_0[ \times \R$}. 
\eq
Towards this end, we point out that the definition of entropy admissible solution~\eqref{e:entropy} is equivalent to requiring for every $k \in \R$ the equations  
 \be \label{e:entropy2}
   \partial_t [\rho- k]^+ + \partial_x \Big[\mathbbm{1}_{\rho \ge k} \big[V(\rho, u) \rho - V(k, u) k \big]\Big] +
   \mathbbm{1}_{\rho \ge k} k\partial_2 V (k, u) \partial_x u \leq 0
\eq
and 
\be \label{e:entropy3}
   \partial_t [\rho- k]^- - \partial_x \Big[\mathbbm{1}_{\rho \leq k} \big[V(\rho, u) \rho - V(k, u) k \big]\Big] -
   \mathbbm{1}_{\rho \leq k} k\partial_2 V (k, u) \partial_x u \leq 0,
\eq 
are both satisfied in the sense of distributions. In the previous expression $[\cdot]^+$ and $[\cdot]^-$ denote the positive and negative part, respectively.  By space-integrating over $\R$ the inequality~\eqref{e:entropy2} with $k=1$ and recalling~\eqref{e:id} and the identity $V(1, u) \equiv0$ we obtain the upper bound in~\eqref{e:bdfrombelow}. By doing the same to~\eqref{e:entropy3} with $k=0$ we get the lower bound.  This handwaving computation can be made rigorous by (for instance) a fairly standard vanishing viscosity approximation argument. \\
{\bf Step 3:} we establish a bound on the total variation of $\rho_{n+1}$. Again, we only provide a formal proof, which can be formally justified through a standard vanishing viscosity argument. We set $\partial_x \rho_{n+1}: = r_{n+1}$, which yields 
\begin{equation}\begin{split} \label{e:erren+1}
   0 & = \partial_t r_{n+1} + \partial_x [V(\rho_{n+1}, u_n) r_{n+1}]+ \partial_x [\partial_x V(\rho_{n+1}, u_n) \rho_{n+1}]
   \\
   & = \partial_t r_{n+1} + \partial_x \Big[ [V(\rho_{n+1}, u_n)+ \partial_1 V \rho_{n+1}] r_{n+1}\Big]
    +
    \partial_x [\partial_2 V \partial_x u_n  \rho_{n+1} ] .
\end{split}
\end{equation}
By recalling that $\partial_x u_n = \rho_n z_n$ we get
\begin{equation*}
\begin{split}
     \partial_x [\partial_2 V \partial_x u_n \rho_{n+1} ] = &  [ \partial_{2 1} V r_{n+1}
     + \partial_{22} V\partial_x u_{n} ] \partial_x u_n  \rho_{n+1} \\ &
    +\partial_2 V \rho_{n+1}
    [\partial_x z_n \rho_n + z_n \partial_x \rho_{n}]+ 
   \partial_2 V  \partial_x u_n r_{n+1} 
\end{split}
\end{equation*}
and by plugging the above formula into~\eqref{e:erren+1}, multiplying times $\mathrm{sign} [r_{n+1}]$ and space integrating we get 
\begin{equation}\label{e:vedirem} \begin{split}
   \frac{d}{dt} &\int_{\R} |r_{n+1}(t, x)| dx  \leq 
   \| \partial_x u_n  \|_{L^\infty} \| \rho_{n+1} \|_{L^\infty}
   \left(  \| \partial_{21 } V \|_{C^0}   \int_{\R} |r_{n+1}(t, x)| dx +   \| \partial_{22 } V \|_{C^0}      \int_{\R} |\partial_x u_{n}(t, x)| dx  \right) \\ &
 +  \| \partial_{2 } V \|_{C^0}  \| \rho_{n+1}  \|_{L^\infty}  \left[
       \| \partial_x z_n \|_{L^\infty}  \| \rho_{n}  \|_{L^1}
    +  \| z_{n}  \|_{L^\infty} \int_{\R} |\partial_x \rho_n(t, x)| dx \right] \\ &+
    \| \partial_{2 } V \|_{C^0}
    \| \partial_x u_n  \|_{L^\infty}
       \int_{\R} |r_{n+1}(t, x)| dx  \\
   & \stackrel{\eqref{e:hpitera1},\eqref{e:hpitera2},\eqref{e:hpitera3}}{\leq} 
    C (V, \|z_0 \|_{L^\infty}) \int_{\R} |r_{n+1}(t, x)| dx + 
    C (V, \| z_0 \|_{L^\infty}, \| \psi_0 \|_{L^\infty}, \| \rho_0 \|_{L^1}) + C (V,    \| z_0 \|_{L^\infty}) M_0. 
\end{split}
\end{equation}
By applying the Comparison Theorem for ODEs we conclude that 
\begin{equation} \label{e:tvrhon}
\begin{split}
     \sup_{t \in [0, \tau_0]} \int_{\R} |r_{n+1}(t, x)| dx & \leq 
    [\mathrm{TotVar} \, \rho_0 ] \exp(\tilde C \tau_0) + [M_0 + 1] [ \exp (\tilde C \tau_0)- 1],  
\end{split}
\end{equation}
for some $\tilde C=  C (V, \| z_0 \|_{L^\infty}, \| \psi_0 \|_{L^\infty}, \| \rho_0 \|_{L^1})$. We now set  
\be \label{e:cienne}
    M_0: = 4 \mathrm{TotVar} \rho_0 + 4 
\eq
and choose $\tau_0$ in such a way that 
\be \label{e:tau0}
    [ \exp (\tilde C \tau_0)- 1] \leq 2 \tau_0 \leq \sfrac{1}{2}  
\eq
and by using~\eqref{e:tvrhon}
we arrive at 
\be \label{e:tvestimate}
        \mathrm{Tot Var} \, \rho_{n+1} (t, \cdot) 
        \leq [\sfrac{1}{2} + 1] \mathrm{TotVar} \rho_0 + \sfrac{1}{2}[ M_0+1] \leq M_0 \quad \text{for $t \in [0, \tau_0]$. }
\eq
{\bf Step 4:} we define $u_{n+1}$. Note that we cannot exactly rely on the classical method of characteristics owing to the low regularity of $\rho_{n+1}$ and henceforth of the vector field $V(\rho_{n+1}, u_n)$.  To circumvent this obstruction, we consider the Cauchy problem obtained by coupling the equation 
\be \label{e:stambecco}
              \partial_t [v_{n+1}] + \partial_x [ V(\rho_{n+1}, u_n)v_{n+1}] =0 
\eq
with the initial datum $v_{n+1} (0, \cdot)= \psi_0 \rho_0$. Given~\eqref{e:kruenne}, existence and uniqueness results for~\eqref{e:stambecco} are available under very weak regularity assumptions on $V(\rho_{n+1}, u_n)$ and $\rho_{n+1}$, see~\cite[\S4]{Panov2008}. In the following for technical reasons we rely on~\cite[Theorem 2.5]{ACFS} applied with $p= \rho_{n+1}$, $b = V(\rho_{n+1}, u_n)$. We conclude that there is a unique solution of the Cauchy problem satisfying $|v_{n+1}| \leq \| \psi_0 \|_{L^\infty} \rho_{n+1}$. Note by setting $\psi_{n+1}: = \sfrac{v_{n+1}}{\rho_{n+1}}$ if $\rho_{n+1}>0$ and $\psi_{n+1} =0$ otherwise we obtain
\be \label{e:arancione}
        v_{n+1} = \psi_{n+1} \rho_{n+1}, \qquad \| \psi_{n+1} \|_{L^\infty} \leq \| \psi_0 \|_{L^\infty}
\eq
Note furthermore that $\| v_{n+1}(t, \cdot) \|_{L^1}\leq \| \psi_{n+1}(t, \cdot) \|_{L^\infty} \| \rho_{n+1}(t, \cdot) \|_{L^1} \leq \| \psi_0 \|_{L^\infty} \| \rho_0 \|_{L^1}$.

To define $z_{n+1}$, we point out that, owing to the identity $z_0' =\rho_0 \psi_0$ and to the $L^1$ and $L^\infty$ bounds on $\rho_0$ and $\psi_0$, respectively, the function $z_0$ is of bounded total variation and hence the limit $z_\infty: = \lim_{x \to - \infty} z_0(x)$ exists and is finite. We then set 
\be \label{e:azzurro}
    z_{n+1} (t, x) : = z_{\infty} + \int_{- \infty}^x \rho_{n+1}\psi_{n+1} (t, x)  dx,
\eq
which is an $L^\infty(\R_+; W^{1, \infty}(\R))$ function owing to the $L^1$ and the $L^\infty$ bounds on $\rho_{n+1}$ and $\psi_{n+1}$, and furtheremore satisfies the initial condition $z_{n+1} = z_0 \psi_0$ due to~\eqref{e:idgarz}. 
Note that, owing to~\eqref{e:stambecco}, we have the identity 
\be \label{e:rosso}
    \partial_t z_{n+1} + V(\rho_{n+1}, u_{n}) \partial_x z_{n+1} =0 \quad \text{a.e. $(t, x) \in \R_+ \times \R$}.
\eq
We now want to establish the estimate 
\be \label{e:verde2}
      |z_{n+1}(t, x)| \leq \| z_0 \|_{L^\infty} \; \text{for every $(t, x) \in \R_+ \times \R$}. 
\eq
If $V(\rho_{n+1}, u_{n})$ were a regular vector field, we could apply the classical method of characteristic and~\eqref{e:verde2} 
would straigthforwardly follow from~\eqref{e:rosso}. Owing to the low regularity of $V(\rho_{n+1}, u_{n})$, the proof of~\eqref{e:verde2} is slightly more involved, and we now provide its details.  

We fix a test function $\varphi$ and use $\varphi z_{n+1}$ as a test function in the definition of distributional solution of~\eqref{e:kruenne}. We conclude that $v_{n+1}= \rho_{n+1} z_{n+1}$ satisfies the very same continuity equation~\eqref{e:stambecco}, now coupled with the initial datum $\rho_0 z_0$. We apply again~\cite[Theorem 2.5]{ACFS} and conclude that  the Cauchy problem obtained by coupling~\eqref{e:stambecco} with the initial datum $\rho_0 z_0$ admits a unique solution satisfying $|v_{n+1}| \leq M \rho_{n+1}$ for some $M>0$, and that one can take $M= \| z_0  \|_{L^\infty}$. By uniqueness, this solution must coincide with $\rho_{n+1} z_{n+1}$, and this implies that 
\be \label{e:verde}
      |z_{n+1}(t, x)| \leq \| z_0 \|_{L^\infty} \; \text{for $\rho_{n+1} \mathcal L^{2}$ a.e. $(t, x) \in \R_+ \times \R$}. 
\eq
We now show that~\eqref{e:verde} implies~\eqref{e:verde2}. We  fix $t>0$, and point out that owing to~\eqref{e:azzurro} $z_{n+1}(t, \cdot)$ is a Lipschitz continuous function. Let $\Lambda_{n+1}(t)$ denote the open set of $x$-s such that $z_{n+1}(t, x) > \| z_0 \|_{L^\infty}$. We decompose $\Lambda_{n+1}(t)$ as a countable union of disjoint intervals, 
$$
    \Lambda_{n+1}(t) = \bigcup_{k=1}^\infty \mathcal I_k.
$$
We now fix $k \in \mathbb N$ and assume by contradiction that $\mathcal I_k \neq \emptyset$. If $\mathcal I_k = \R$, then owing to~\eqref{e:verde} we have $\rho_{n+1}(t,  \cdot) =0$ a.e on $\R$, and owing to~\eqref{e:azzurro} this yields $z_{n+1}(t, \cdot) = z_\infty$. Since $z_\infty$ is the asymptotic state of $z_0$, then $|z_\infty | \leq \| z_0 \|_{L^\infty}$ and hence this contradicts the definition of $\Lambda_{n+1}(t) $. If $\mathcal I_k \neq \R$, then $\mathcal  I_k = ] a_k, b_k [$ and at least one between $a_k$ and $b_k$ is finite. Just to fix the ideas, let us assume that $a_k \in \R$, then we observe that owing to~\eqref{e:verde} we must have $\rho_{n+1}=0$ a.e. on $\Lambda_{n+1}(t)$. We conclude that for every $x \in \mathcal I_k$ we have 
\begin{equation*} \begin{split}
     |z_{n+1} (t, x) | & = \left| z_{n+1} (t, a_k) + \int_{a_k}^x v_{n+1} (t, x) dx \right| \leq 
     |z_{n+1} (t, a_k) | + \int_{a_k}^x |\psi_{n+1}| \underbrace{\rho_{n+1}(t, x) }_{=0}dx\\
    &  = |z_{n+1} (t, a_k) | \leq \| z_0 \|_{L^\infty}, 
\end{split}
\end{equation*}
where in the last equality we have used that $a_k \notin \Lambda_{n+1}(t)$. This contradicts the definition of $\Lambda_{n+1}(t)$ an hence establishes~\eqref{e:verde2}. 

We now set 
\be \label{e:uennepiu}
   u_{n+1} (t, x) : = u_\infty +  \int_{-\infty}^x z_{n+1} \rho_{n+1}(t, x) dx, 
\eq
where $u_\infty$ is the limit at $-\infty$ of function $u_0$, which is of bounded total variation owing to the identity $u_0' = \rho_0 z_0$. By following the same argument as before we get\footnote{Note that the bound from below in~\eqref{e:bdperiterare3} is not mentioned in the statement of~\cite[Theorem 2.5]{ACFS}, but it is a straightforward consequence of the proof, given that $u_0 \ge 0$.}
\be \label{e:bdperiterare3}
    0 \leq  u_{n+1}(t, \cdot)  \leq \| u_0 \|_{L^\infty}
\eq 
and also 
\be \label{e:te2}
    \partial_t u_{n+1} + V(\rho_{n+1}, u_{n+1}) \partial_x u_{n+1}=0
     \quad \text{a.e $(t, x) \in \R_+ \times \R$. }
\eq
Note furthermore that the previous analisis yields  
\be \label{e:bdperiterare}
   \partial_x u_{n+1} =  \rho_{n+1} z_{n+1}
\eq
and 
\be \label{e:bdperiterare2}
    \partial_x z_{n+1} = \rho_{n+1} \psi_{n+1}, \quad \| \psi_{n+1} \|_{L^\infty} \leq \| \psi_0 \|_{L^\infty}. 
  \eq
{\bf Step 5:} by combining~\eqref{e:rhon+1elle1},~\eqref{e:bdfrombelow} and \eqref{e:tvestimate} we get that $\rho_{n+1}$ satisfies~\eqref{e:hpitera1}, by combining~\eqref{e:bdperiterare},\eqref{e:bdperiterare3} and~\eqref{e:bdperiterare2} we get that $u_{n+1}$ satisfies~\eqref{e:hpitera2} and~\eqref{e:hpitera3}, which implies that we can iterate the construction for every $n \in \mathbb N$. 
\subsection{Passage to the limit}\label{sss:lpass}
We  fix the same $\tau_0$ as in~\eqref{e:tau0}. We recall the bounds~\eqref{e:bdfrombelow} and~\eqref{e:tvestimate} and control the time derivative by using the equation at the first line of~\eqref{e:kruenne} combined with the Volpert Chain Rule and the bound on $\partial_x u_n$ obtained by combining~\eqref{e:verde2} and~\eqref{e:bdperiterare}. We apply the Helly-Fr\'echet-Kolmogorov Compactness Theorem and we conclude that there is $\{ \rho_{n_k} \}$ such that 
\be \label{e:magenta}
    \rho_{n_k} \to \rho \quad \text{strongly in $L^1 ([0, \tau_0] \times \R$)}, 
\eq
for some limit function $\rho$ such that $\mathrm{TotVar} \, \rho(t, \cdot) \leq M_0$ for a.e. $t \in [0, \tau_0].$  Note that the above convergence occurs in $L^1$ and not only in $L^1_{\mathrm{loc}}$ because owing to~\eqref{e:tortora} for every $t \in [0, \tau_0]$ the function $\rho_{n+1}(t, \cdot)$ vanishes outside the interval $[- \tilde R, \tilde R]$ for some $\tilde R$ depending on $R$ and $\tau_0$ but independent of $n$.  

Next, we recall the bound~\eqref{e:bdperiterare}, the identity $\partial_x z_{n+1} = \rho_{n+1} \psi_{n+1}$ (which follows from~\eqref{e:azzurro}) and the bound in~\eqref{e:arancione}. We also use the equation~\eqref{e:rosso} to deduce a uniform bound on the time derivative. We apply the Arzel\`a Ascoli Theorem and conclude that there is subsequence (which to simplify the notation we do not relabel) such that $\{  z_{n_k} \}$ converges to some limit function $z$ uniformly on the compact subsets of $[0, \tau_0] \times \R$. Note furthermore that $\partial_x z_{n_k}$ an $\partial_t z_{n_k}$ converge to $\partial_x z$ and $\partial_t z$ weakly$^\ast$ in $L^\infty ([0, \tau_0] \times \R)$ and that  we have the identity $\partial_x z= \rho \psi$, where $\psi$ is an accumulation point of the sequence $\{ \psi_n \}_{n \in \mathbb N}$ in the weak$^\ast$ topology and as such satisfies the bound in~\eqref{e:verde3}. 

By combining~\eqref{e:magenta} with the uniform convergence of $\{ z_{n_k} \}$ we can pass to the limit in the identity~\eqref{e:uennepiu} and conclude that the sequence $\{ u_{n_k} \}$ converges uniformly on $[0, \tau_0] \times \R$ to the limit function 
\be \label{e:u}
   u (t, x) : = u_\infty +  \int_{-\infty}^x z \rho(t, x) dx.
\eq
By passing to the limit in~\eqref{e:bdperiterare3} we obtain the first two inequalities in~\eqref{e:verde3}. Assume for a moment that we have shown that the sequence $\{ u_{n_k -1} \}$ converges to the very same limit $u$ given by~\eqref{e:u}\footnote{By compactness, $\{ u_{n_k -1} \}$ converges up to subsequences to some limit function $v$, but in principle $u$ and $v$ could be different.} Then we can pass to the limit in the distributional formulation (and in the entropy inequality) for the equation at the first line of~\eqref{e:kruenne} and obtain an entropy admissible solution of the equation at the first line of~\eqref{e:GARZ}. By passing to the limit in the identity~\eqref{e:te2} we also obtain the equation at the second line of~\eqref{e:GARZ}. Note furtheremore that the identity $\partial_x u = z \rho$ and the bound $\| z \|_{L^\infty} \leq \| z_0 \|_{L^\infty}$ directly follow from~\eqref{e:u} and from~\eqref{e:verde2} and the uniform convergence of $\{ z_{n_k} \}$. 

The above argument establishes existence of an admissible solution of the Cauchy problem obtained by coupling~\eqref{e:GARZ} with~\eqref{e:idpose} defined on the time interval $[0, \tau_0]$. To define a global in time solution, it suffices to point out that the value of the constant $\tilde C$ in~\eqref{e:tvrhon} and henceforth the value of $\tau_0$ in~\eqref{e:tau0} \emph{does not} depend on the total variation of the initial datum. By using the bounds~\eqref{e:verde} and $\| \rho(t, \cdot) \|_{L^1} \leq \| \rho_0 \|_{L^1}$ we conclude that the constants $\tilde C$ and $\tau_0$ only depend on quantities that are preserved by the admissible solution and hence that we can iterate the construction and in this way establish global in time existence. 

To conclude the proof of the existence part of Theorem~\ref{t:wpl} we are thus left to show  that the sequence $\{ u_{n_k -1} \}$ converges to the same limit as $\{ u_{n_k } \}$. Towards this end, we introduce the functional 
$$
    \Phi_n (t) : =\|[\rho_n -\rho_{n+1}](t,\cdot)\|_{L^1} + \| [\partial_x u_{n} - \partial_x u_{n-1}](t,\cdot) \|_{L^1}  
$$
and introduce the following result. 
\begin{lemma}
\label{l:phienne}
We have  
\be\label{e:phienne2}
   \lim_{n  \to + \infty} \Phi_n (t ) \to 0 \quad \text{for every $t \ge 0$}.
\eq
\end{lemma}
We pospone the proof of Lemma~\ref{l:phienne} to the next paragraph and we now show that~\eqref{e:phienne2} implies the sequence $\{ u_{n_k -1} \}$ converges to the same limit as $\{ u_{n_k } \}$ on $[0, \tau_0] \times \R$. We
recall that, owing to~\eqref{e:uennepiu}, $\lim_{x \to - \infty}u_n(t, x) = u_\infty$ for every $n$ and every $t \in [0, \tau_0]$, which in particular implies 
\be 
\label{e:linftyl1}
    \| u_{n} - u_{n-1} \|_{C^0} \leq \| \partial_x u_{n} - \partial_x u_{n-1} \|_{L^1}.
\eq 
By combining~\eqref{e:linftyl1} with~\eqref{e:phienne2} we then conclude that $\lim_{n \to + \infty} \| u_{n}(t, \cdot) - u_{n-1} (t, \cdot)\|_{C^0}=0$ for every $t \ge 0$.
\subsection{Proof of Lemma~\ref{l:phienne}}\label{ss:pphienne}
We apply the stability result proven in \cite[Theorem 2.6]{MRS} (see also Remark 2.8 therein about entropy solutions) for conservation laws with space-time dependent fluxes $P=P(t,x,\rho), \ Q= Q(t,x,\rho)$. In our setting the fluxes take the form
\[
(t,x,\rho) \mapsto P(t,x,\rho) \doteq \rho_{n+1} V(\rho_{n+1}, u_n(t,x)), \qquad (t,x,\rho) \mapsto Q(t,x,\rho) \doteq \rho_{n+1} V(\rho_n, u_{n-1}(t,x)).
\]
and the stability estimate in \cite{MRS} implies that for every $t\ge 0$ we have 
\begin{equation}\label{e:stab_cl}
\begin{split}
\int_{\R}& |\rho_{n+1}(t,x)-\rho_n(t,x)|  dx \le ~  \underbrace{\int_{\R} |\rho_{n+1}(0,x)-\rho_n(0,x)| dx}_{=0} \\
& + C(V) \int_{0}^{t} (\|z_n (s)\|_{L^\infty} + \|z_{n-1}(s)\|_{L^\infty} ) \int_\R |\rho_{n+1}(s,x) - \rho_n(s,x)| dx ds \\
& + C(V) \underbrace{\int_{0}^{t} \int_\R \left(|\partial_x u_n| |u_n-u_{n-1}| + |\partial_x u_n - \partial_x u_{n-1}| \right) dx ds}_{: = J}  + C(V) M_0  \int_{0}^{t} \underbrace{\|(u_n(s, \cdot)-u_{n-1}(s, \cdot))\|_{L^\infty}}_{\leq \|[\partial_x u_{n} - \partial_x u_{n-1}](s, \cdot) \|_{L^1} \; \text{by \eqref{e:linftyl1}}} ds
\end{split}
\end{equation}
where $M_0$ is the same as in~\eqref{e:hpitera1}. We now control the term $J$ as follows: 
\be \label{e:notheta}
\begin{split}
    J &\leq \int_{0}^{t} \| [u_n- u_{n-1} ](s, \cdot) \|_{L^\infty} \underbrace{\int_\R |\partial_x u_n |(s, x) dx}_{\leq \| z_0 \|_{L^\infty} \| \rho_0 \|_{L^1}} ds + \int_{0}^{t} \|[\partial_x u_{n} - \partial_x u_{n-1}](s, \cdot) \|_{L^1}ds \\
   &  \stackrel{\eqref{e:linftyl1}}{\leq} C(\| z_0 \|_{L^\infty} ,\| \rho_0 \|_{L^1})
    \int_{0}^{t} \|[\partial_x u_{n} - \partial_x u_{n-1}](s, \cdot) \|_{L^1}ds,
\end{split}
\eq
and conclude that 
\be \label{e:pezzo1}
  \int_{\R} |\rho_{n+1}(t,x)-\rho_n(t,x)|  dx \leq
  C(V, M_0, \| z_0 \|_{L^\infty} ,\| \rho_0 \|_{L^1})  \int_{0}^{t} \Phi_{n}(s)ds 
\eq
We now point out that  
\be \label{e:pezzo2} \begin{split}
  &
\|[\partial_x u_n- \partial_xu_{n-1}](t, \cdot) \|_{L^1}
  =   \| [ \rho_n z_n  - \rho_{n-1} z_{n-1}](t, \cdot) \|_{L^1}
   \\ & \leq  \| \rho_n [ z_n  -  z_{n-1}] (t, \cdot) \|_{L^1} + \| z_{n-1} \|_{L^\infty}
    \|  [\rho_n   - \rho_{n-1} ](t, \cdot)  \|_{L^1}\\
& \stackrel{\eqref{e:pezzo1}}\leq \| \rho_n [ z_n  -  z_{n-1}]  \|_{L^1} +
   C(V, M_0, \| z_0 \|_{L^\infty} ,\| \rho_0 \|_{L^1})  \int_{0}^{t} \Phi_{n-1}(s)ds 
\end{split}
\eq
and that
\[
   \partial_t [\rho_n z_n] + \partial_x [V(\rho_n, u_{n-1}) \rho_n z_n ] =0, \quad 
    \partial_t [\rho_{n-1} z_{n-1}] + \partial_x [V(\rho_{n-1}, u_{n-2}) \rho_{n-1} z_{n-1} ] =0, 
\]
which implies (by using the equation for $\rho_n$ and $\rho_{n-1}$ and the Volpert Chain Rule)
\begin{equation*}
\begin{split}
   0 &= \partial_t \big[\rho_n [z_n-z_{n-1}]\big] + \partial_t \big [z_{n-1} [\rho_n - \rho_{n-1}] \big] + 
   \partial_x \big[V(\rho_n, u_{n-1}) \rho_n [z_n - z_{n-1}] \big] \\& \quad + 
    \partial_x \big[z_{n-1} [ V(\rho_n, u_{n-1}) \rho_n - V(\rho_{n-1}, u_{n-2}) \rho_{n-1}] \big]\\
  & =  \partial_t \big[\rho_n [z_n-z_{n-1}]\big] + \partial_t z_{n-1}[\rho_n - \rho_{n-1}]+ 
        z_{n-1} \partial_t [\rho_n - \rho_{n-1}]
    +   \partial_x \big[V(\rho_n, u_{n-1}) \rho_n [z_n - z_{n-1}] \big]\\ & \qquad + 
      \partial_x z_{n-1} [ V(\rho_n, u_{n-1}) \rho_n - V(\rho_{n-1}, u_{n-2}) \rho_{n-1}] +
     z_{n-1} \partial_x \big[ V(\rho_n, u_{n-1}) \rho_n - V(\rho_{n-1}, u_{n-2}) \rho_{n-1}] \\
   & \stackrel{\eqref{e:kruenne}}{=} 
    \partial_t \big[\rho_n [z_n-z_{n-1}]\big] + \partial_t z_{n-1}[\rho_n - \rho_{n-1}]
    +   \partial_x \big[V(\rho_n, u_{n-1}) \rho_n [z_n - z_{n-1}] \big]\\ & \qquad + 
      \partial_x z_{n-1} [ V(\rho_n, u_{n-1}) \rho_n - V(\rho_{n-1}, u_{n-2}) \rho_{n-1}] \\
     & \stackrel{\eqref{e:rosso}}{=}
      \partial_t \big[\rho_n [z_n-z_{n-1}]\big] - V(\rho_{n-1}, u_{n-2}) \partial_x z_{n-1} [\rho_n - \rho_{n-1}]
        +   \partial_x \big[V(\rho_n, u_{n-1}) \rho_n [z_n - z_{n-1}] \big]  \\ & \quad + 
      \partial_x z_{n-1} [ V(\rho_n, u_{n-1}) \rho_n - V(\rho_{n-1}, u_{n-2}) \rho_{n-1}] .
\end{split}
\end{equation*}
By recalling the equation~\eqref{e:kruenne} for $\rho_n$ and the Volpert Chain Rule we infer that the Lipschitz continuous function $z_n-z_{n-1}$ satisfies 
\begin{equation*}\begin{split}
     0 & =  \rho_n \big[ \partial_t  [z_n-z_{n-1}] + V(\rho_n, u_{n-1})  \partial_x  [z_n-z_{n-1}] \big] - V(\rho_{n-1}, u_{n-2}) \partial_x z_{n-1} [\rho_n - \rho_{n-1}]
        \\ & + 
      \partial_x z_{n-1} [ V(\rho_n, u_{n-1}) \rho_n - V(\rho_{n-1}, u_{n-2}) \rho_{n-1}], 
      \end{split}
\end{equation*}
which by multiplying times  $\mathrm{sign}(z_n - z_{n-1})$ yields 
\begin{equation*}\begin{split}
     0 & =  \rho_n \big[ \partial_t  |z_n-z_{n-1}| + V(\rho_n, u_{n-1})  \partial_x  |z_n-z_{n-1}|  \big] - \mathrm{sign}(z_n - z_{n-1}) V(\rho_{n-1}, u_{n-2}) \partial_x z_{n-1} [\rho_n - \rho_{n-1}]
        \\ & + 
     \mathrm{sign}(z_n - z_{n-1}) \partial_x z_{n-1} [ V(\rho_n, u_{n-1}) \rho_n - V(\rho_{n-1}, u_{n-2}) \rho_{n-1}] =0.
      \end{split}
\end{equation*}
By using again the equation for $\rho_n$ we then arrive at 
\begin{equation*}\begin{split}
     0 & = \partial_t \big[\rho_n |z_n-z_{n-1}| \big] -  \mathrm{sign}(z_n - z_{n-1}) V(\rho_{n-1}, u_{n-2}) \partial_x z_{n-1} [\rho_n - \rho_{n-1}]
        +   \partial_x \big[V(\rho_n, u_{n-1}) \rho_n |z_n - z_{n-1}| \big]  \\ & \quad + 
       \mathrm{sign}(z_n - z_{n-1}) \partial_x z_{n-1} [ V(\rho_n, u_{n-1}) \rho_n - V(\rho_{n-1}, u_{n-2}) \rho_{n-1}]
      \end{split}
\end{equation*}
and by integrating the above inequality in space and time we arrive at
\begin{equation}\label{e:stab_tr}
\begin{split}
   \int_{\R} & \rho_n |z_n-z_{n-1}|(t, \cdot)  dx \leq
   C(V) \| \partial_x z_{n-1} \|_{L^\infty} 
   \left( \int_0^t \int_\R |\rho_n - \rho_{n-1}| dx ds 
  +  \int_0^t \int_\R \rho_n  | u_{n-1} - u_{n-2} |  dx ds  \right) \\
  & \stackrel{\eqref{e:bdfrombelow},\eqref{e:arancione}}{\leq}
    C(V, \| \psi_0 \|_{L^\infty} )
    \int_0^t \int_\R |\rho_n - \rho_{n-1}|dx ds \\ &
  +   C(V, \| \psi_0 \|_{L^\infty} ) \int_0^t \| [ u_{n-1}(s, \cdot) - u_{n-2}(s, \cdot)) \|_{L^\infty}\underbrace{\int_\R \rho_n(s, x) dx}_{\leq \| \rho_0 \|_{L^1}} ds \\
    & \stackrel{\eqref{e:linftyl1}}{\leq}
       C(V, \| \psi_0 \|_{L^\infty}, \| \rho_0 \|_{L^1} ) \int_0^t \Phi_{n-1} (s) ds  
\end{split}
\end{equation}
By recalling~\eqref{e:pezzo1} and~\eqref{e:pezzo2} we eventually conclude that
$$
   \Phi_n(t) \leq
  C(V, M_0, \| z_0 \|_{L^\infty} ,\| \rho_0 \|_{L^1})  \int_{0}^{t} \Phi_{n}(s)ds+
  C(V, \| \psi_0 \|_{L^\infty})  \int_{0}^{t} \Phi_{n-1}(s)ds,  
$$
which owing to the Gr\"onwall Lemma implies that for every $0 \leq \tau \leq t$ we have 

\be \label{e:growiter}
\begin{split}
    \Phi_n (t) & \leq \Phi_n (\tau) \exp\big( C(V, M_0, \| z_0 \|_{L^\infty} ,\| \rho_0 \|_{L^1}) [t - \tau]\big)
    \\
     & +
    C(V, M_0, \| z_0 \|_{L^\infty} ,\| \rho_0 \|_{L^1})  [\exp \big( C(V, M_0, \| z_0 \|_{L^\infty} ,\| \rho_0 \|_{L^1}) [t - \tau] 
     \big)  -1 ] \sup_{s \in [\tau, t] }\Phi_{n-1}(s). 
\end{split}
\eq
We now argue iteratively: first, we choose $\tau= 0$ and $\delta>0$ in such a way that 
$$
[\exp \big( C(V, M_0, \| z_0 \|_{L^\infty} ,\| \rho_0 \|_{L^1} \delta)-1] \leq \sfrac{1}{2}. 
$$
Owing to~\eqref{e:growiter} and to the equality $\Phi_n(0)=0$ this implies $\sup_{t \in [0, \delta]} \Phi_n (t) \leq \sfrac{1}{2} \sup_{t \in [0, \delta] }\Phi_{n-1}(t)$ and hence, by induction on $n$, $\sup_{t \in [0, \delta]} \Phi_n (t) \leq (\sfrac{1}{2})^n \sup_{t \in [0, \delta]} \Phi_0 (t)$. Next, we plug this inequality into~\eqref{e:growiter} applied with $\tau=\delta$, and get 
$$
    \sup_{t \in [\delta, 2 \delta]} \Phi_n (t) \leq \sfrac{3}{2} \left(\sfrac{1}{2} \right)^n \sup_{t \in [0, \delta]} \Phi_0 (t) + \sfrac{1}{2} 
   \sup_{t \in [\delta, 2\delta]} \Phi_{n- 1} (t),
$$
which by induction yields $  \sup_{t \in [\delta, 2 \delta]} \Phi_n (t)  \leq (n+1) (\sfrac{1}{2})^n [ \sfrac{3}{2} \sup_{t \in [0, 2 \delta]}
\Phi_0 (t)]$. By iterating the above argument we arrive at~\eqref{e:phienne2}.

\section{Uniqueness and stability}\label{ss:luni2}
We follow the same argument as in the proof of Lemma~\ref{l:phienne} and we provide the details for the sake of completeness. 
We fix $T>0$ and $t\in [0, T]$ and we apply~\cite[Theorem 2.6]{MRS} with 
\[
(t,x,\rho) \mapsto P(t,x,\rho) \doteq \rho V(\rho, u_1(t,x)), \qquad (t,x,\rho) \mapsto Q(t,x,\rho) \doteq \rho V(\rho, u_2(t,x)),
\]
which yields 
\begin{equation}\label{e:stab_cl2}
\begin{split}
\int_{\R}|\rho_1(t,x)-\rho_2(t,x)|  dx \le &~  \int_{\R}|\rho_1(0,x)-\rho_2(0,x)| dx \\
& + C(V) \int_{0}^{t} (\|z_1(s)\|_{L^\infty} + \|z_2(s)\|_{L^\infty} ) \int_\R |\rho_1(s,x) - \rho_2(s,x)| dx ds \\
& + C(V, \| \rho_1 \|_{L^\infty}, \| \rho_2 \|_{L^\infty}) \underbrace{\int_{0}^{t} \int_\R \left(|\partial_x u_1| |u_1-u_2| + |\partial_x u_1 - \partial_x u_2| \right)dx ds}_{:=L} \\
& + C(V, \| \rho_1 \|_{L^\infty}, \| \rho_2 \|_{L^\infty}) M_0  \int_{0}^{t} \|(u_1(s, \cdot)-u_2(s, \cdot))\|_{L^\infty} ds
\end{split}
\end{equation}
where $M_0$ is such that $ \mathrm{TotVar} \, \rho_1(s, \cdot) \le M_0$ and $ \mathrm{TotVar} \, \rho_2(s, \cdot) \le M_0$ for every $t \in [0, T]$. By using the identity $\partial_x u_1 = \rho_1 z_1$ and the estimate $\| \rho_1 (t, \cdot) \|_{L^1} \leq \| \rho_0 \|_{L^1}$ we control the term $L$ as follows: 
\be \label{e:notheta2}
    L \leq \| \rho_0 \|_{L^1} \| z_1 \|_{L^\infty} \int_{0}^{t} \| [u_1- u_2 ](s, \cdot) \|_{L^\infty} 
     + \int_{0}^{t} \|[\partial_x u_1 - \partial_x u_2](s, \cdot) \|_{L^1}ds. 
\eq
We now want to show that 
\be \label{e:limitiuguali}
      \lim_{x \to - \infty }  u_1 (t, x) =  \lim_{x \to - \infty } u_1 (0, x), 
     \quad
      \lim_{x \to - \infty }  u_2 (t, x) =  \lim_{x \to - \infty } u_2 (0, x)  \quad \text{for every $t>0$}
\eq
Note that all the above limit exist and are finite because  $u_1(t, \cdot)$ and $u_2 (t, \cdot)$ are all functions of bounded variation owing to the identities $\partial_x u_i(t, \cdot) = \rho_i z_i (t, \cdot)$ and to the $L^1$ and $L^\infty$ bounds on $\rho_i$ and $z_i$, respectively.  By contradiction, assume that there is $t>0$ such that one of the equalities in~\eqref{e:limitiuguali} fails, for instance $\lim_{x \to - \infty }  u_1 (t, x)  \neq  \lim_{x \to - \infty }u_1 (0, x) (x)$. This yield the existence of $t>0$ and $d>0$ such that 
\be \label{e:tortora2}
    \int_{-R}^{-R+1} |u_1(t, x) - u_1(0, x)| dx \ge d \quad \text{for every $R$ sufficiently large}.
\eq
On the other hand, by using the equation at the first line of~\eqref{e:GARZ} we have 
\begin{equation*}
\begin{split}
         \int_{-R}^{-R+1} |u_1(t, x) - u_1(0, x)| dx  & =  \int_{-R}^{-R+1} \left| \int_0^t \partial_\tau u_1(\tau, x) d \tau \right| dx   \stackrel{\eqref{e:GARZ}}{\leq}
              \int_0^t \int_{-R}^{-R+1}  |V (\rho_1, u_1) \partial_x u_1 (\tau, x) | d \tau dx \\
              & \leq C(V)  \int_0^t \int_{-R}^{-R+1}  |\partial_x u_1 (\tau, x) | d \tau dx \to 0 \; \text{as $R \to + \infty$.}
\end{split}
\end{equation*}
To establish the convergence at the last line of the above expression we have used the Lebesgue Dominated Convergence Theorem combined with the identity $\partial_x u_1(t, \cdot) = \rho_1 z_1 (t, \cdot)$, which dictates (owing to the bounds on $\rho_1$ and $z_1$) that $\partial_x u_1$ is a bounded and summable function. The above convergence contradicts~\eqref{e:tortora2} and establishes~\eqref{e:limitiuguali}. By using~\eqref{e:limitiuguali} we then get
\be \label{e:indaco}
    \|(u_1(s, \cdot)-u_2(s, \cdot))\|_{L^\infty}
   \leq  |u_{1\infty}- u_{2 \infty}| + \|\partial_x u_1(s, \cdot)-\partial_x u_2(s, \cdot)\|_{L^1},
\eq
where $u_{i \infty}$ denotes the asymptotic state of $u_i(0, \cdot)$ as $x \to - \infty$. 
We now point out that  
\be \label{e:cremisi}
\|\partial_x u_1- \partial_xu_2 \|_{L^1}
 =   \|  \rho_1 z_1  - \rho_2 z_2 \|_{L^1}
    \leq  \| \rho_1 [ z_1  -  z_2]  \|_{L^1} + \| z_2 \|_{L^\infty}
    \|  \rho_1   - \rho_2   \|_{L^1}
\eq
so we are actually left to control $\| \rho_1 [ z_1  -  z_2]  \|_{L^1}$.
By combining the equations at the first and second line of~\eqref{e:GARZ} with the identity $\partial_x u_i = \rho_i z_i$ we get
\be \label{e:celeste}
   \partial_t [\rho_1 z_1] + \partial_x [V(\rho_1, u_1) \rho_1 z_1 ] =0, \quad 
    \partial_t [\rho_2 z_2] + \partial_x [V(\rho_2, u_2) \rho_2 z_2 ] =0.
\eq
We also have 
\be \label{e:bianco}
    \rho_1 [\partial_t z_1 + V(\rho_1, u_1) \partial_x z_1 ] =0, 
    \quad 
     \rho_2 [\partial_t z_2 + V(\rho_2, u_2) \partial_x z_2 ] =0,
\eq
This yields  
\begin{equation} \label{e:beige}
\begin{split}
   0 &\stackrel{\eqref{e:celeste}}{=} \partial_t \big[\rho_1 [z_1-z_2]\big] + \partial_t \big [z_2 [\rho_1 - \rho_2] \big] + 
   \partial_x \big[V(\rho_1, w_1) \rho_1 [z_1 - z_2] \big] + 
    \partial_x \big[z_2 [ V(\rho_1, u_1) \rho_1 - V(\rho_2, u_2) \rho_2] \big]\\
  & =  \partial_t \big[\rho_1 [z_1-z_2]\big]  + \partial_t z_2 [\rho_1 - \rho_2] + z_2 \partial_t [\rho_1-\rho_2]
    +  \partial_x \big[V(\rho_1, u_1) \rho_1 [z_1 - z_2] \big] \\ & \qquad + 
      \partial_x z_2 [ V(\rho_1, u_1) \rho_1 - V(\rho_2, u_2) \rho_2] +
     z_2 \partial_x \big[ V(\rho_1, u_1) \rho_1 - V(\rho_2, u_2) \rho_2] \big] \\
   & \stackrel{\eqref{e:GARZ}} = 
      \partial_t \big[\rho_1 [z_1-z_2]\big] + \partial_t z_2 [\rho_1 -  \rho_2]
      +  \partial_x \big[V(\rho_1, u_1) \rho_1 [z_1 - z_2] \big]  + 
      \partial_x z_2 [ V(\rho_1, u_1) \rho_1 - V(\rho_2, u_2) \rho_2].
\end{split}
\end{equation}
Note furthermore that 
\begin{equation*}
\begin{split}
          \partial_t z_2 [\rho_1 -   \rho_2]  &=
          \rho_1 \partial_t [z_2 - z_1] + \rho_1 \partial_t z_1 - \partial_t z_2 \rho_2
            \stackrel{\eqref{e:bianco}}{=} 
           \rho_1 \partial_t [z_2 - z_1] - \partial_x z_1 V(\rho_1, u_1) \rho_1+
         \partial_x z_2 V(\rho_2, u_2) \rho_2 \\
         & =  \rho_1 \partial_t [z_2 - z_1] + V(\rho_2, u_2) \rho_2 \partial_x [z_2 - z_1]
         + \partial_x z_1 [ V(\rho_2, u_2) \rho_2 - V(\rho_1, u_1) \rho_1]
\end{split}
\end{equation*}
and by plugging the above equation into~\eqref{e:beige}, multiplying the result by $\mathrm{sign}(z_1 - z_2)$ and arguing as in the proof of Lemma~\ref{l:phienne} we 
arrive at 
\be \label{e:smeraldo}
\begin{split}
 0 &=   \partial_t \big[\rho_1 |z_1-z_2| \big]+  \rho_1 \partial_t |z_2 - z_1| + V(\rho_2, u_2) \rho_2 \partial_x |z_2 - z_1|
         + \mathrm{sign}(z_1 - z_2)\partial_x z_1 [ V(\rho_2, u_2) \rho_2 - V(\rho_1, u_1) \rho_1]
      \\ & \qquad +  \partial_x \big[V(\rho_1, u_1) \rho_1 |z_1 - z_2| \big]  
       + \mathrm{sign}(z_1 - z_2)      \partial_x z_2 [ V(\rho_1, u_1) \rho_1 - V(\rho_2, u_2) \rho_2].
\end{split}
\eq 
By using the Integration by Parts Formula combined with the equation for $\rho_1$ (or more rigorously by combining a suitable approximation argument with the definition of distributional solution for $\rho_1$) we get  
\[ 
    \begin{split}
    \int_0^t \int_{\R} \rho_1 \partial_t |z_2 - z_1| (s, x) dx ds & =
    \int_{\R} \rho_1 |z_2 - z_1| (t, x) - 
   \int_{\R} \rho_1 (0, x) |z_2 - z_1| (0, x) \\ & - \int_0^t \int_{\R} V(\rho_1, u_1)  \rho_1 \partial_x |z_2 - z_1| (s, x) dx ds.
   \end{split}
\]
By space and time integrating~\eqref{e:smeraldo} and using the above identity  we arrive at 
\begin{equation}\label{e:aprile}\begin{split}
    \int_{\R} & \rho_1 |z_1 - z_2|  (t, x) dx \leq
     \int_{\R}  \rho_1 |z_1 - z_2|  (0, x) dx \\
    & + 
    \int_0^t \int_{\R} |V(\rho_1, u_1) \rho_1- V(\rho_2, u_2) \rho_2| [| \partial_x z_2|+  |\partial_x z_1|] (s, x) dx ds \\
   & \leq   \int_{\R}  \rho_1 |z_1 - z_2|  (0, x) dx   + 
     C(V, \| \partial_x z_1 \|_{L^\infty},  \| \partial_x z_2 \|_{L^\infty} )
      \Bigg[ \int_0^t \| [u_1 - u_2 ](s, \cdot) \|_{L^\infty} \underbrace{ \int_{\R} \rho_1 (s, x) dx}_{\leq \| \rho_0 \|_{L^1}} ds  \\ & \quad
      +
      \int_0^t \| [\rho_1 - \rho_2 ](s, \cdot) \|_{L^1} ds\Bigg] \\
   & \stackrel{\eqref{e:indaco}}{\leq}   \int_{\R}  \rho_1 |z_1 - z_2|  (0, x) dx  
       + C(V, \| \partial_x z_1 \|_{L^\infty},  \| \partial_x z_2 \|_{L^\infty}, \| \rho_1(0, \cdot) \|_{L^1} )
      \Bigg[  t | u_{1 \infty} - u_{2 \infty}|   \\
     & \qquad +  \int_0^t \| [\partial_x u_1 - \partial_x u_2 ](s, \cdot) \|_{L^1}  ds + 
      \int_0^t \| [\rho_1 - \rho_2 ](s, \cdot) \|_{L^1} ds.
    \Bigg]
\end{split}
\end{equation}
Next, we point out that 
\be \label{e:maggio}
      \int_{\R}  \rho_1 |z_1 - z_2|  (0, x) dx \stackrel{0 \leq \rho_{1}(0, \cdot) \leq 1}{\leq} \| \partial_x u_1 (0, \cdot) - 
      \partial_x u_2 (0, \cdot) \|_{L^1} + \| z_2(0, \cdot) \|_{L^\infty} \| [\rho_1- \rho_2 ] (0, \cdot) \|_{L^1}
\eq
and we set 
\[
\Phi(t)=   \|[\rho_1-\rho_2](t,\cdot)\|_{L^1} + \| [\partial_x u_1 - \partial_x u_2](t,\cdot) \|_{L^1} .
\]
By combining~\eqref{e:stab_cl2},~\eqref{e:notheta2}~\eqref{e:indaco},~\eqref{e:cremisi},~\eqref{e:aprile} and \eqref{e:maggio}, and using $\|\partial_x u_i\|_{L^\infty}\le \| \rho_i \|_{L^\infty} \|z_i\|_{L^\infty}$, we have 
\[
\Phi(t) \le \hat C \Big[ \Phi(0) + |u_{1\infty} - u_{2 \infty}| T+  \int_0^t \Phi(s) ds  \Big], \quad \text{for every $t \in [0, T]$}, 
\]
where the constant $\hat C$ depends on $V$, the total variation of $\rho_1, \rho_2$ and $\|z_1\|_{L^\infty}, \|  \rho_1(0, \cdot) \|_{L^1}, \|z_2\|_{L^\infty}, \|\partial_x z_2\|_{L^\infty}$, $\| \rho_1 \|_{L^\infty}$ and $\| \rho_2 \|_{L^\infty}$.
Owing to the Gronwall Lemma this yields 
\[
    \Phi(t) \le e^{\hat C t}   [\Phi(0) + T |u_{1\infty} - u_{2 \infty}|],
\]
that is~\eqref{e:stability}. The above estimate yields in particular the identities $u_1=u_2$ and $\rho_1=\rho_2$ if the solutions have the same initial data.
\section*{Acknowledgments}
Both authors are members of the GNAMPA group of INDAM and are supported by the project PRIN 2022 PNRR C53D23008420001 (PI Roberta Bianchini). LVS is also supported by PRIN 2022YXWSLR (PI Paolo Antonelli) and by the CNR project STRIVE (DIT.AD022.207 - STRIVE (FOE 2022)). Both PRIN projects are financed by the European Union-Next Generation EU.  \bibliographystyle{plain}
\bibliography{garz2}
\end{document}